\documentstyle[11pt]{article}

\newtheorem{th}{Theorem}[section]

\newtheorem{Theorem}[th]{Theorem}

\newtheorem{Lemma}[th]{Lemma}

\newtheorem{Corollary}[th]{Corollary}
\newtheorem{Definition}[th]{Definition}

\newcommand{\qed}{\hbox{\vrule\vbox to10pt{\hrule width
 4pt\vfill\hrule}\vrule}\smallskip}
\newenvironment{proof}{\noindent {\bf Proof.\/}}{\hfill\qed}

\newcommand{\bbr}{{\bf R}}
\newcommand{\reals}{{\bf R}}

\newcommand{\integers}{{\bf Z}}
\newcommand{\rationals}{{\bf Q}}

\newcommand{\pcirc}{\kern .7pt {\scriptstyle \circ} \kern 1pt}

\newcommand{\eqref}[1]{(\ref{#1})}

\newcommand{\pol}{{\rm Pol\,}}
\newcommand{\Td}{{\rm Td\,}}

\newcommand\lie[1]{{\bf #1}}
\newcommand{\onehalf}{\frac{1}{2}}

\begin{document}

\title{A limit of toric symplectic forms that has no periodic Hamiltonians}
\author{Jean-Claude HAUSMANN
\footnote{University of Geneva, Switzerland; hausmann@ibm.unige.ch}
\and Allen KNUTSON
\footnote{Brandeis University, Waltham, Massachusetts;
 allenk@alumni.caltech.edu}
\footnote{Both authors thank the Fonds National Suisse de la
Recherche Scientifique for its support.}}
\date{\today}

\maketitle
\begin{abstract}
We calculate the Riemann-Roch number of some of the pentagon spaces defined
in \cite{Kl,KM,HK1}.
Using this, we show that while the regular pentagon space is
diffeomorphic to a toric variety, even symplectomorphic to one under
arbitrarily small perturbations of its symplectic structure,
it does not admit a symplectic circle action.
In particular, within the cohomology classes of symplectic structures,
the subset admitting a circle action is not closed.
\end{abstract}

\section{The result}

In \cite{Kl,HK1,KM} are introduced {\em polygon spaces}, which are
a family of
symplectic manifolds often possessing many interesting circle actions.
In the case of pentagons with distinct edge-lengths, they are actually
(symplectomorphic to) toric varieties. So it is rather a surprise that
the regular pentagon space, although {\em diffeomorphic} to the nearby
nearly-regular pentagon spaces, admits no symplectic circle action at all.

In what follows we calculate the Riemann-Roch number of the
regular pentagon space, $6$, and its Euler characteristic, $7$.
This will already be enough to show it isn't
a toric variety. To prove the stronger statement we invoke results of
Karshon \cite{Ka} on arbitrary symplectic $4$-manifolds with circle actions.

\section{Nearly-regular pentagon spaces and toric structures}

We remind the reader of the definition of pentagon spaces and why
they are frequently symplectomorphic to toric varieties \cite{HK1,KM}.
In what follows our cohomology groups will always be with real coefficients.

\subsection{Pentagon spaces.}
Given a list $(a_1,\ldots,a_5)$ of positive real numbers, define
the {\em pentagon space} $\pol(a_1,\ldots,a_5)$ as the space of 5-step
closed piecewise-linear paths in $\reals^3$ starting at the origin,
with step-lengths $a_1,\ldots,a_5$, considered up to $SO(3)$-rotation.
As such, it is a subquotient of the {\em path space}
$\prod_{i=1}^5 S^2_{a_i}$ of 5-step paths in $\reals^3$ starting at the
origin, not necessarily closing up nor quotiented by $SO(3)$.

Regarding $S^2_{a_i}$ as a symplectic manifold with area $a_i$, this
construction can be interpreted as producing the
pentagon space as a symplectic quotient of the path space
by $SO(3)$. Consequently, the pentagon space is also symplectic when it
is smooth; this happens exactly if no two of the edges put together
have the same total length as the other three put together.

\subsection{Nearly-regular pentagon spaces.}
\begin{Definition}
We define a {\em nearly-regular} pentagon space $\pol(a_1,\ldots,a_5)$
to be one in which the sum of any two edge-lengths $a_i+a_j$ is {\em less}
than the sum of the other three.
\end{Definition}

In particular these spaces are smooth, and of course include
the regular pentagon space $\pol(1,1,1,1,1)$.

As we will see below, any nearly-regular pentagon space
is diffeomorphic to $(S^2\times S^2) \# \overline{CP^2} \# \overline{CP^2} \# \overline{CP^2}$;
so when we now restrict attention to nearly-regular pentagon spaces,
we can view this as studying different symplectic structures on
the same manifold.

\subsection{Toric structures.}
Take now the generic case that $a_1\neq a_2$, $a_4\neq a_5$.

In \cite{Kl,KM} are defined the {\em bending flows} on such a polygon space.
Label the vertices A,B,C,D,E, and draw the diagonals AC and AD,
breaking the pentagon into three triangles ABC, ACD, and ADE.
Since these diagonals cannot be zero length (by the genericity assumption
above), each one gives a circle action: break off the triangles ABC and ADE
and reattach them rotated by independent angles.
In fact these circle actions on this symplectic manifold are Hamiltonian,
with the Hamiltonians given by the lengths of the diagonals $(d_1,d_2)$.

%
\begin{figure}[Hht]
  \begin{center}
    \leavevmode
    \input{pent.pstex_t}
  \end{center}
\end{figure}
%
In \cite[p. 192]{HK1} we calculated the image of the map $(d_1,d_2)$ in $\reals^2$;
in the nearly-regular case we always get a rectangle missing three corners.
(See figure \ref{fig:momentpoly} below.)
With this and Delzant's reconstruction theorem it is immediate that the
space is diffeomorphic to
$(S^2\times S^2) \# \overline{CP^2} \# \overline{CP^2} \# \overline{CP^2}$.
In particular the $7$ vertices says that the Euler characteristic is $7$.

\begin{figure}[Hht]
    \leavevmode
    \input{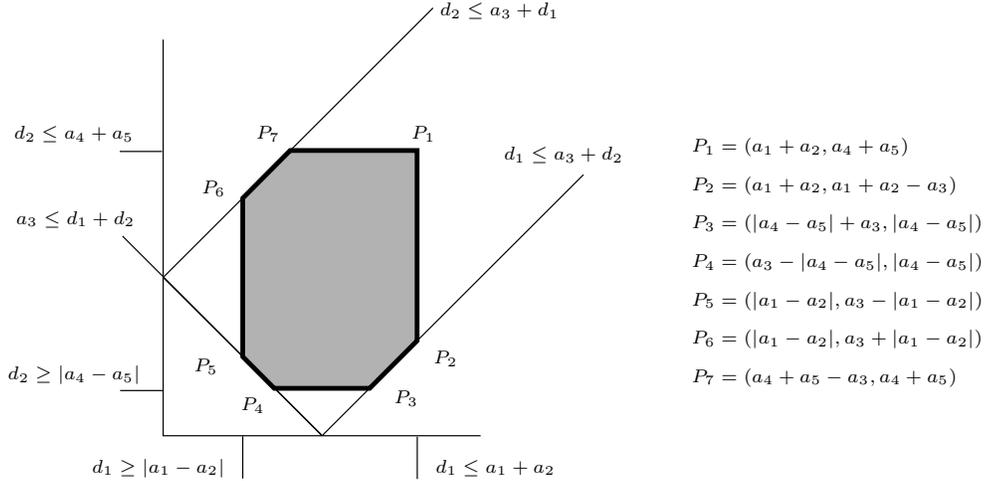}
\caption{The moment polytope of $\pol(a_1,\ldots,a_5)$.}
    \label{fig:momentpoly}
\end{figure}

Actually one can do slightly better than this abstract diffeomorphism.
Fix $(M,\omega) := Pol(1,1,1,1,1)$. 
Then for $a:=(a_1,...,a_5)$ the edge-lengths of a
nearly-regular pentagon space there is a canonical (up to isotopy)
diffeomorphism $h: M \to Pol(a_1,...,a_5)$ (\cite[Prop. 2.2]{HK2}).
Call $\omega_a\in H^2(M)$ the (well-defined) class $\omega_a:=h^*\omega$.

\subsection{Riemann-Roch numbers of nearly-regular pentagon spaces.}
Any compact symplectic manifold $(M,\omega)$ has a canonical ``Todd class''
$\Td M$,
constructed by endowing the manifold with an auxiliary almost
complex structure. (In fact pentagon spaces are actually complex \cite{Kl}.)
The {\em Riemann-Roch number} $RR(M,\omega)$ is defined as
$\int_M \exp([\omega]) \Td M$ and only depends on the cohomology class
$[\omega]$ of the symplectic form. 

Note that we have not assumed integrality of $\omega$ 
for this definition; rather, this is a polynomial function on
$H^2(M)$ which merely happens to take on integral values at
integral $\omega$. In the case of pentagon spaces this exactly corresponds
to the side-lengths being integral \cite[corollary 6.6, p. 304]{HK2}.

\begin{Lemma}\label{intrr}
As in the last section, 
fix $M$ the regular pentagon space, 
but change the symplectic form to $\omega_a$, where
$\pol(a_1,\ldots,a_5)$ is a toric, integral, nearly-regular pentagon space.
Then the Riemann-Roch number of $(M,\omega_a)$ is
$$ \onehalf(a_1+a_2+a_3+a_4+a_5)^2 - 2 (a_1^2+a_2^2+a_3^2+a_4^2+a_5^2)
+ \onehalf(a_1+a_2+a_3+a_4+a_5) + 1.$$
\end{Lemma}

\begin{proof}
If $M$ is a toric manifold with moment map $\mu:M \to \lie{t}^*$,
$\omega$ integral, and some (therefore every) fixed point of $M$ maps
to a lattice point in $\lie{t}^*$ under $\mu$, then the Riemann-Roch
number of $M$ is equal to the number of lattice points in the
image of $\mu$ \cite{Fu}.

In our case the moment map is $\mu=(d_1,d_2)$.
Assume for convenience that $a_1>a_2$ and $a_4>a_5$ (they may be
switched if necessary without affecting the Riemann-Roch number).
Then by Pick's theorem \cite[p. 113]{Fu} the
number of lattice points in figure \ref{fig:momentpoly} is the area plus
one-half the number of lattice points on the boundary plus one.
The area is
$$4 a_2 a_5
- \onehalf(a_2+a_4+a_5-a_1-a_3)^2
- \onehalf(a_2+a_3+a_5-a_1-a_4)^2
- \onehalf(a_1+a_2+a_5-a_3-a_4)^2 $$
and the boundary has
$$ 2(2 a_2+1) + 2(a_1-a_2)-1 + a_3+a_4+a_5-a_1-a_2-1 = a_1+a_2+a_3+a_4+a_5 $$
lattice points. Expanding and collecting terms we get the result
(which is, as it should be, symmetric in $a_1\ldots a_5$).
\end{proof}

\begin{Corollary}\label{5goRR}
The formula from lemma \ref{intrr} 
holds for any nearly-regular pentagon space; 
in particular $RR(\pol(1,1,1,1,1)) = 6$.
\end{Corollary}

\begin{proof}
Let $C$ be any open cone in $\reals^n$ of dimension $n$,
i.e. an open subset closed under multiplication by positive reals.
(For us, the relevant $\reals^n$ will be $H^2(M)$, and $C$ the space of
classes of symplectic forms from nearly-regular pentagon spaces.)

We first observe that a polynomial $p$ on $\reals^n$ is determined by its
values on $C \cap \integers^n$.
Proof: $p$ is dominated by its highest-order terms, $p_{top}$. 
If $p_{top}$ vanishes on $C \cap \integers^n$, it vanishes
on $C \cap \rationals^n$ by homogeneity, 
therefore on all of $C$ by continuity, therefore everywhere 
by analyticity. 

Next we find a basis of $H^2(M)$ with which to apply this.
Recall that we constructed the pentagon space by reduction from
the $5$-step path space $(S^2)^5$ by $SO(3)$.
Since $H^2(BSO(3))$ is trivial (recall we're using
real coefficients), by Kirwan's theorem the linear
map from $H^2_{SO(3)}((S^2)^5)$ to $H^2(M)$ is
surjective. By dimension count it is an isomorphism (the moment polytope
lets us determine that $b_2(M)$ is also $5$).
So the obvious basis of $H^2((S^2)^5)$ -- which is what we were
using to write down the vectors $(a_1,\ldots,a_5)$ --
gives us a basis for $H^2(M)$.

Since the Riemann-Roch number of $(M,\omega)$ is a polynomial function on
$H^2(M)$, and has been determined on the subset of integral, 
positive, distinct $(a_1,\ldots,a_5)$ in lemma \ref{intrr},
the formula in lemma \ref{intrr} in fact calculates the 
Riemann-Roch number for any nearly-regular pentagon space.
\end{proof}

It is worth mentioning a corollary: the symplectic volume (which is 
always the leading term in the Riemann-Roch number) 
of a nearly-regular pentagon space is
$\onehalf(a_1+a_2+a_3+a_4+a_5)^2 - 2 (a_1^2+a_2^2+a_3^2+a_4^2+a_5^2)$.

\section{The regular pentagon space has no periodic Hamiltonians} \label{}
Let $P:=\pol (1,1,1,1,1)$. So far we know its Riemann-Roch number is $6$
and its Euler characteristic is $7$.

\begin{Lemma}
$P$ is not a toric manifold.
\end{Lemma}

\begin{proof}
Suppose that $P$ admits a Hamiltonian action of a 2-torus $T$.
Let $\Pi$ be the moment polytope in ${\rm Lie\,}(T)^*$.
As $P$ is an integral symplectic manifold, one may translate $\Pi$ so
that its vertices are lattice points.
Since $P$'s Euler characteristic is $7$, $\Pi$ has $7$ vertices, and
thus at least $7$ lattice points.
This contradicts the fact that $RR(P)=6$,
since the Riemann-Roch number of a toric manifold is the number of
lattice points in its moment polytope.
\end{proof}

\begin{Theorem}\label{5gono}
$P$ admits no periodic hamiltonian.
\end{Theorem}

\begin{proof}
Suppose that $f:P\to\bbr$ is a moment map for the
Hamiltonian action of a circle on $P$.

One may suppose that the minimum value of $f$ is $0$.
As $P$ is an integral symplectic manifold, the critical values of
$f$ are then integers (see \cite[Lemma 2.5]{Ka} and its proof).
Observe that there are more than two critical values. Otherwise,
$P$ would be $CP^2$ or an $S^2$-bundle over $S^2$ (\cite[2.4.2]{Au}),
neither of which fit with $b_2(P)=5$.
Also, at least one of the extremal level sets
(say the minimal one) is not an isolated point.
Otherwise, $f$ would have only isolated
critical points and then be a toric manifold by \cite[Theorem 5.1]{Ka},
which contradicts the lemma above.

We now consider the Duistermaat-Heckman function $DH$ for $f$
which is a piecewise affine map from $\bbr$ to $\bbr_{\geq0}$.
We shall use the following facts:
\begin{enumerate}
\item $DH(0)$ is the area of the sphere $f^{-1}(0)$ which
is an integer since $\omega$ is integral.
\item The slope of $DH$
at $0$ is the Euler class of the normal bundle
of $f^{-1}(0)$ in $P$ (\cite[Lemma 2.12]{Ka}).
\item the slope of $DH$  changes exactly at each critical value of $f$ and the
graph of $DH$ is concave (\cite[Lemmas 2.12 and 2.19]{Ka}.
\item \label{last} the number of integral points
enclosed by the graph of $DH$ is \\ $RR(P)=6$.
\end{enumerate}
One checks that the only possible $DH$ function must then satisfy
$$ DH(0)=1, \quad  DH(1)=2 \  \hbox{ and }  DH(y)=0 \  \hbox{ for }
y\geq 2.$$
Therefore, the maximum level $f^{-1}(2)$ would be an isolated point.
But the slope of $DH$ at an extremal isolated point should be $\leq 1$
in absolute value (\cite[Lemma 2.12]{Ka}) which is not the case here.
\end{proof}

\vskip .3 truecm\goodbreak

\vskip .5 truecm\small
\noindent \parbox[t]{6 truecm}{Jean-Claude HAUSMANN\\
Math\'ematiques-Universit\'e\\ B.P. 240, \\
 CH-1211 Gen\`eve
 24, Suisse\\ hausmann@math.unige.ch} \ \hfill \hfill \
\parbox[t]{5 truecm}{Allen KNUTSON \\
 Department of Mathematics\\
 Brandeis University\\
 Waltham, MA 02254-9110 USA\\
 allenk@alumni.caltech.edu}

\end{document}